\documentclass[12pt,a4paper]{article}

\usepackage{epsf,epsfig,amsfonts,amsgen,amsmath,amstext,amsbsy,amsopn,amsthm
}
\usepackage{amsmath,times,mathptmx}
\usepackage{amsfonts,amsthm,amssymb}
\usepackage{amsfonts}
\usepackage{graphics}
\usepackage{latexsym,bm}
\usepackage{amsfonts,amsthm,amssymb,bbding}
\usepackage{indentfirst}
\usepackage{graphicx}
\usepackage{color}
\usepackage[colorlinks=true,anchorcolor=blue,filecolor=blue,linkcolor=blue,urlcolor=blue,citecolor=blue]{hyperref}
\usepackage{float}
\usepackage{tikz}
\usepackage{verbatim}
\usepackage{mathrsfs}

\evensidemargin=\oddsidemargin\topmargin=-1.5cm

\newtheorem{theorem}{Theorem}[section]
\newtheorem{prob}{Problem}[section]

\newtheorem{lemma}{Lemma}[section]
\newtheorem{cor}{Corollary}[section]

\newtheorem{claim}{Claim}[section]

\addtocounter{section}{0}

\begin{document}
\title{Spectral Radius of Graphs with Size Constraints: Resolving a Conjecture of Guiduli
\footnote{Supported by the National Natural Science Foundation of China (No. \!12171066).}
}
\author{{ Rui Li$^a$}, \quad  {Anyao Wang$^b$}, \quad { Mingqing Zhai$^c$}\thanks{Corresponding author. E-mail addresses: lirui@hhu.edu.cn (R. Li),
 way@ahnu.edu.cn (A. Wang), mqzhai@njust.edu.cn (M. Zhai).}\\
{\footnotesize $^a$ School of Mathematics, Hohai University, Nanjing 211100, Jiangsu, P.R. China}\\
{\footnotesize $^b$ School of Mathematics and Statistics, Anhui Normal University,
         Wuhu 241002, Anhui, P.R. China}\\
{\footnotesize $^c$ School of Mathematics and Statistics, Nanjing University of Science and Technology}\\
{\footnotesize Nanjing 210094, Jiangsu, P.R. China}}

\date{}
\maketitle

\begin{abstract}

We resolve a problem posed by Guiduli (1996) on the spectral radius of graphs satisfying the Hereditarily Bounded Property $P_{t,r}$, which requires that every subgraph $H$ with $|V(H)| \geq t$ satisfies $|E(H)| \leq t|V(H)| + r$. For an $n$-vertex graph $G$ satisfying $P_{t,r}$, where $t > 0$ and $r \geq -\binom{\lfloor t+1 \rfloor}{2}$, we prove that the spectral radius $\rho(G)$ is bounded above by  
$\rho(G) \leq c(s,t) + \sqrt{\lfloor t \rfloor n}$,
where $s = \binom{\lfloor t \rfloor + 1}{2} + r$, thus affirmatively answering Guiduli's conjecture.  

Furthermore, we present a complete characterization of the extremal graphs that achieve this bound. 
These graphs are constructed as the join graph $K_{\lfloor t \rfloor} \nabla F$, where $F$ is either $K_3 \cup (n - \lfloor t \rfloor - 3)K_1$ or a forest consisting solely of star structures. The specific structure of such forests is meticulously characterized. 
 Central to our analysis is the introduction of a novel potential function $\eta(F) = e(F) + (\lfloor t \rfloor - t)|V(F)|$, which quantifies the structural "positivity" of subgraphs. By combining edge-shifting operations with spectral radius maximization principles, we establish sharp bounds on $\eta^+(G)$, the cumulative positivity of $G$. 
Our results contribute to the understanding of spectral extremal problems under edge-density constraints and provide a framework for analyzing similar hereditary properties. 

\medskip


\noindent {\bf AMS\,\,\,Classification:} 05C50; 05C35

\end{abstract}

\section{Introduction}
All graphs considered in this context are simple, undirected and free from isolated vertices.
We use $|G|$ and $e(G)$ to denote the numbers of vertices and edges in a graph
$G$, respectively. For a subset $S$ of $V(G)$, let $G[S]$ be the subgraph induced by $S$.
Let $N_U(v)$ denote the set of neighbors of vertex $v$ within the subset $U$ of vertices,
formally defined as
$ N_U(v) = \{ u \mid u \in U \text{ and } uv \in E(G) \}$.
In particular, when $U=V(G)$, we simply write $N(v)$ for $N_{V(G)}(v)$.
Denoted by $d_U(v)$ the cardinality of the set $N_U(v)$.
Moreover, let $e(U)$ denote the size of $G[U]$ and
$e(U,W)$ be the number of edges connecting $U$ and $W$ in $G$.
Let $kG$ be the disjoint union of $k$ copies of $G$, and
$G_1\nabla G_2$ be the graph obtained by joining each vertex of $G_1$ with each vertex of $G_2$.

Let $A(G)$ denote the adjacency matrix of a graph $G$,
and let $\rho(G)$ represent the spectral radius of $A(G)$.
A graph $G$ is said to have the {\it Hereditarily Bounded Property}
$P_{t,r}$ if for each subgraph $H \subseteq G$ with $|V(H)| \geq t$,
the size of $H$ satisfies $ |E(H)| \leq t \cdot |V(H)| + r$, where $t>0$.
In fact, many classes of graphs possess the Hereditarily Bounded Property,
including planar graphs,
which have been extensively studied in the context of spectral extremal problems
(see \cite{boots91, cao93, Cvetkovic90, mohar10, zha00, lin21, zhai}).
Research into the spectral radius of graphs with similar hereditary properties has also flourished.
For instance, constraints on matching numbers have been explored by Feng et al. \cite{feng07},
while constraints on cycle lengths were investigated by Gao et al. \cite{gao19}.
These studies contribute to a broader understanding of how the structural properties
of graphs influence their spectral characteristics.

This concept was first introduced by Guiduli in his Ph.D. thesis \cite{Guiduli96},
wherein he proved the following theorem.
\begin{theorem}\label{t0}
Let $t\in \mathbf{N}$ and $r\geq -\binom{ t+1 }{2}$.
If $G$ is a graph on $n$ vertices with the property $P_{t,r}$, then
$\rho(G)\leq (t-1)/2+\sqrt{ t(t+1)+2r}+\sqrt{ tn}$.
\end{theorem}

For ease of notation, Guiduli defined $s = \binom{t+1}{2}+r$ in his thesis.
When $t$ and $ r$ are not restricted to integers,
this result implies that $\rho(G) \leq c(s,t) + \sqrt{\lceil t \rceil n}$,
where $c(s,t)$ represents a function of $s$ and $t$.
Guiduli also posed the following problem.

\begin{prob}
  Let $t,r\in \mathbf{R}$, $t>0$ and $r\geq -\binom{\lfloor t+1 \rfloor}{2}$.
  Is it true that $G$ has the property $P_{t,r}$, then $\rho(G)\leq c(s,t)+\sqrt{ t  n}$?
\end{prob}

Furthermore, Guiduli conjectured a potentially stronger statement:
\[\rho(G) \leq c(s,t) + \sqrt{\lfloor t \rfloor n}.\]

In this paper, we affirmatively resolve the conjecture posed by Guiduli.
Specifically, we establish the following theorem:
\begin{theorem}\label{t2}
Let $t,r\in \mathbf{R}$, $t>0$ and $r\geq -\binom{\lfloor t+1 \rfloor}{2}$.
If $G$ is an $n$-vertex graph with the property $P_{t,r}$ , then
$\rho(G)\leq c(s,t)+\sqrt{\lfloor t\rfloor n}$.
\end{theorem}

Finally, we consider the structure of $G$
when it attains the maximum spectral radius over all
$n$-vertex graphs satisfying the property $P_{t,r}$.
The following theorem provides a precise characterization of such graphs.

\begin{theorem}\label{t3}
Let $G$ be a graph that attains the maximum spectral radius over all
$n$-vertex graphs satisfying the property $P_{t,r}$. Then $G$ is the join of two graphs:
one is $K_{\lfloor t\rfloor}$,
another is either $K_3\cup(n-\lfloor t\rfloor-3)K_1$ or a forest consisting of stars.
\end{theorem}

In Section 4,
we provide a detailed characterization of the forest structure in Theorem \ref{t3},
classifying the stars into maximal, big, small, and minimal stars based on their sizes.
Therefore, all the extremal graphs are completely characterized.

\section{Stability theorem}
Let $\mathcal{G}$ denote the family of $n$-vertex graphs that satisfy property $P_{t,r}$.
For any $G \in \mathcal{G}$, we assume for simplicity that the graph $G$
considered in this paper is connected. If $G$ were not connected,
we could identify a component $C$ within $G$ such that the induced subgraph $G[C]$
also satisfies property $P_{t,r}$ and has the same spectral radius as $G$, i.e.,
$\rho(G[C]) = \rho(G)$.
Therefore, it suffices to consider only the connected component $G[C]$ for our analysis.
Furthermore, we may assume that $n$ is sufficiently large.
If $n$ were a finite constant, we could set $c(s,t) = n$.
In this case, the spectral radius would be bounded by
$\rho(G)\leq\rho(K_n)=n - 1 \leq c(s,t) + \sqrt{\lfloor t \rfloor n}$.
Hence, Theorem \ref{t2} would hold trivially.

Thus, for the remainder of our discussion,
we focus on connected graphs $G \in \mathcal{G}$ with a sufficiently large order $n$.
Given that $G\in \mathcal{G}$, it implies that $t\geq 1$. To see why,
suppose for contradiction that $0 < t < 1$. For sufficiently large $n$,
$ e(G) \leq tn + r < n - 1$.
However, since $G$ is connected, $e(G)\geq n-1$.
This leads to a contradiction, hence $t \geq 1 $.

By the well-known Perron-Frobenius theorem,
there exists a positive unit eigenvector $X $
corresponding to the spectral radius $\rho(G) $.
Let $X = (x_1, x_2, \ldots, x_n)^T$
and let $u^* \in V(G)$ such that $x_{u^*} = \max_{u \in V(G)} x_u.$
For simplicity, we will write $\rho = \rho(G)$.

For convenience, we define
$s=\big\lceil\frac{t(\lfloor t\rfloor+1)+r}{\lfloor t\rfloor+1-t}\big\rceil+1$ in this paper.
This definition of $s$ differs from that given in Guiduli's thesis
but similarly depends on $t$ and $r$. It is clear that $s > t + 1$
as $r\geq -\binom{\lfloor t \rfloor+1}{2}$.
To prove our main result,
we introduce a stability theorem that plays a crucial role in our analysis.

\begin{theorem}\label{t1}
Let $G$ be a graph of sufficiently large order $n$. Let $X$ be
a non-negative eigenvector corresponding to $\rho(G)$, and $x_{u^{*}}=\max_{u\in V(G)}x_u$.
If $G$ has the property $P_{t,r}$ and $\rho(G)\geq \sqrt{\lfloor t\rfloor(n-\lfloor t\rfloor)}$,
then there exists a vertex subset $V$ of $G$ with $|V|=\lfloor t\rfloor$
such that $x_v\geq (1-\frac{1}{50s^2})x_{u^{*}}$
and $d_G(v)\geq (1-\frac{1}{25s^2})n$ for every $v\in V$.
\end{theorem}

Next, we will prove Theorem \ref{t1}.
We'll define the set of vertices with large entries in $X$ relative to $x_{u^*}$,
provide an upper bound for its cardinality, and then determine its exact value.
This method was first introduced by Cioab\u{a} et al.
to address a conjecture by Nikiforov \cite{Nikiforov10},
which comprises two parts:
a spectral version of Erd\H{o}s-S\'{o}s conjecture \cite{tait23}
and the spectral even cycle problem \cite{tait24}.
It was further developed by Byrne et al. \cite{Byrne}.

Let $V_1=\{v\in V(G):\ x_v\geq (5s)^{-4} x_{u^*}\}$, $V_2=\{v\in V(G):\ x_v
\geq (5s)^{-5} x_{u^*}\}$ and $V=\{v\in V(G):\ x_v\geq (5s)^{-1} x_{u^*}\}$.
Obviously, $V\subseteq V_1\subseteq V_2$.

\begin{claim}\label{c3}
  $|V_1|\leq (5s)^{-6}n$ and $|V_2|\leq (5s)^{-5}n$.
\end{claim}

\begin{proof}
Let $v \in V_1$. We start with the inequality
\[ \rho x_v \geq \sqrt{\lfloor t\rfloor(n-\lfloor t\rfloor)}(5s)^{-4}x_{u^*}. \]
Moreover, given that $n$ is sufficiently large, it follows that
$\rho x_v\geq2s(5s)^{6} x_{u^*}.$
Additionally, we have
$\rho x_v = \sum_{u \in N(v)} x_u \leq |N(v)|x_{u^*},$
where $N(v)$ denotes the set of neighbors of $v$ in $G$. Therefore,
\[ \sum_{v \in V_1} \rho x_v \leq \sum_{v \in V_1} |N(v)|x_{u^*}
\leq \sum_{v \in V(G)} |N(v)|x_{u^*} = 2e(G)x_{u^*}. \]

Since $G \in \mathcal{G}$ and $s>t+1$,
we know that $e(G) \leq tn + r \leq sn$.
Combining this with the previous inequalities, we obtain
\[ 2s(5s)^{6}|V_1|\leq\sum_{v \in V_1} \rho x_v \leq 2e(G)\leq 2sn. \]
Consequently, $ |V_1| \leq(5s)^{-6}n$, as desired.
Similarly, $|V_2| \leq (5s)^{-5}n.$
\end{proof}

For any vertex $u \in V(G)$, let $N_2(u)$ denote the set of all vertices
that are at a distance of two from $u$ in the graph $G$.
We can categorize the vertices around $u$
into several sets based on their distances from $u$
and their memberships in $V_1$ and $V_2$.
Specifically, with respect to $V_1$, we denote
\begin{itemize}
  \item $A_1 = V_1 \cap N(u)$,  $B_1 = N(u) \setminus A_1$,
  $A_2 = V_1 \cap N_2(u)$,  $B_2 = N_2(u) \setminus A_2$.
\end{itemize}

Similarly, concerning the set $V_2$, we have
\begin{itemize}
  \item $C_1 = V_2 \cap N(u)$,  $D_1 = N(u) \setminus C_1$, $C_2 =
  V_2 \cap N_2(u)$,  $D_2 = N_2(u) \setminus C_2$.
 \end{itemize}

\begin{claim}\label{c4}
For every $u\in V(G)$, we have
\begin{eqnarray}\label{a1}
\lfloor t\rfloor(n-\lfloor t\rfloor)x_u\leq |N(u)|x_u + \Big(\frac{2sn}{(5s)^6}
+\frac{2sn}{(5s)^4}\Big)x_{u^*}
+\sum_{v\in B_1, w\in N_{A_1\cup A_2}(v)} x_w.
\end{eqnarray}
\end{claim}
\begin{proof}
Given that $\rho(G) \geq \sqrt{\lfloor t\rfloor(n-\lfloor t\rfloor)}$ and an arbitrary
vertex $u \in V(G)$, we have
\begin{equation}\label{e1-0}
\lfloor t\rfloor(n-\lfloor t\rfloor) x_{u} \leq \rho^2 x_{u}=\sum_{v\in N(u)} \rho x_v
= |N(u)|x_u + \sum_{v\in N(u), w\in N(v)\setminus\{u\}} x_w.
\end{equation}

Note that $N(u) = A_1 \cup B_1$. Then $ N(v)\setminus\{u\}=
N_{A_1\cup A_2 }(v)\cup N_{B_1\cup B_2}(v)$ for each $v \in N(u)$.
By the definition of $V_1$,
it follows that $x_w < (5s)^{-4}x_{u^*}$ for every $w \in B_1 \cup B_2$.
Therefore, we can establish an upper bound for the vertices in $A_1$ as follows:
\begin{align}\label{e1-1}
\sum_{v\in A_1, w\in N(v)\setminus \{u\}} x_w
&= \sum_{v\in A_1, w\in N_{A_1\cup A_2}(v)} x_{u^*} +
\sum_{v\in A_1, w\in N_{B_1\cup B_2}(v)} (5s)^{-4}x_{u^*} \nonumber\\
&\leq \big(2e(A_1)+e(A_1,A_2)\big)x_{u^*} + e(A_1,B_1\cup B_2)(5s)^{-4}x_{u^*}.
\end{align}
Similarly, for vertices in $B_1$, we obtain
\begin{equation}\label{e1-2}
\sum_{v\in B_1, w\in N(v)\setminus\{u\}} x_w \leq
\sum_{v\in B_1, w\in N_{A_1\cup A_2}(v)} x_w + \big(2e(B_1) + e(B_1, B_2)\big)(5s)^{-4}x_{u^*}.
\end{equation}
Since $G$ satisfies the property $P_{t,r}$ and $s>t+1$, we have that
$2e(A_1) + e(A_1, A_2)\leq 2s|V_1|$ and
$e(A_1, B_1\cup B_2) + 2e(B_1) + e(B_1, B_2)\leq 2sn.$
Additionally, by Claim \ref{c3}, we know that $|V_1|\leq(5s)^{-6}n$.
Combining these three inequalities with our previous results (\ref{e1-0}-\ref{e1-2}),
we can derive the desired inequality (\ref{a1}) directly.
\end{proof}

Using an argument analogous to that in Claim \ref{c4}, with respect to $V_2$,
we can find that for every $u\in V(G)$,
\begin{eqnarray}\label{a2}
\lfloor t\rfloor(n-\lfloor t\rfloor)x_u\leq |N(u)|x_u + \frac{4sn}{(5s)^5}
x_{u^*}
+\sum_{v\in D_1, w\in N_{C_1\cup C_2}(v)} x_w.
\end{eqnarray}

\begin{claim}\label{c5}
$|V_1| < (5s)^6$.
\end{claim}

\begin{proof}
We begin by proving that $|N(u)| \geq (5s)^{-5}n$ for each $u \in V_1$.
Suppose, for the sake of contradiction,
that there exists a vertex $u_0 \in V_1$ with $|N(u_0)| < (5s)^{-5}n$.

By inequality (\ref{a2}), we have
\begin{eqnarray}\label{a3}
\lfloor t \rfloor(n - \lfloor t \rfloor)x_{u_0} \leq |N(u_0)|x_{u_0} + \frac{4sn}{(5s)^5}x_{u^*}
+ \sum_{v \in D_1, w \in N_{C_1 \cup C_2}(v)} x_w.
 \end{eqnarray}
Recall that $D_1 \subseteq N(u_0)$ and $C_1 \cup C_2 \subseteq V_2$.
By the property $P_{t,r}$, we know that
\[ e(D_1, C_1 \cup C_2) \leq e(N(u_0) \cup V_2) < s(|N(u_0)| + |V_2|). \]
Given our assumption, $|N(u_0)| < (5s)^{-5}n$, and by Claim \ref{c3}, $|V_2| \leq (5s)^{-5}n$.
Hence,
\[ |N(u_0)|x_{u_0} + \sum_{v \in D_1, w \in N_{C_1 \cup C_2}(v)} x_w \leq \Big(|N(u_0)|
+ e(D_1, C_1 \cup C_2)\Big)x_{u^*} < \frac{(2s+1)n}{(5s)^5}x_{u^*}. \]

Combining this with inequality (\ref{a3}), we get that
$\lfloor t \rfloor(n - \lfloor t \rfloor)x_{u_0} < \frac{6s+1}{(5s)^{5}}nx_{u^*}.$
On the other hand, since $\lfloor t \rfloor \geq 1$ and
$x_{u_0} \geq (5s)^{-4}x_{u^*}$ (because $u_0 \in V_1$), we have
\[ \lfloor t \rfloor(n - \lfloor t \rfloor)x_{u_0}
\geq \frac{7}{10}nx_{u_0} > \frac{7s}{(5s)^{5}}nx_{u^*}, \]
leading to a contradiction.
Therefore, $|N(u)| \geq (5s)^{-5}n$ for each vertex $u \in V_1$.

Summing this inequality over all vertices $u \in V_1$, we obtain
\[ |V_1|(5s)^{-5}n \leq \sum_{u \in V_1}|N(u)| \leq \sum_{u \in V(G)}|N(u)| = 2e(G) < 2sn. \]
Thus, we conclude that $ |V_1|<(5s)^6$.
\end{proof}

\begin{claim}\label{c6}
For each $u\in V_1 $,
we have $|N(u)|\ge\big(\frac{x_u}{ x_{u^*}}-\frac{1}{(5s)^3}\big)n$.
\end{claim}

\begin{proof}
Consider $S$ as a subset of vertices from $B_1$, with each vertex in $S$
having at least $\lfloor t \rfloor$ neighbors within $A_1\cup A_2$.
Let $p = |A_1| + |A_2|$. We aim to prove that the size of $S$ is
bounded by $|S| \leq s \cdot \ell$, where $\ell = \binom{p}{\lfloor t \rfloor}$.

If $p \leq \lfloor t \rfloor - 1$, then it's impossible for any vertex
in $B_1$ to have $\lfloor t \rfloor$ neighbors in $A_1 \cup A_2$,
which implies that $S$ must be empty, and the statement holds trivially.
Therefore, we focus on the case where $p \geq \lfloor t \rfloor$.

Suppose, for contradiction, that $|S| > s \cdot \ell$.
Each vertex in $S$ connects to at least $\lfloor t \rfloor$ vertices in $A_1 \cup A_2$.
Considering all possible subsets of $A_1 \cup A_2$ of size $\lfloor t \rfloor$,
there are $\ell= \binom{p}{\lfloor t \rfloor}$ such subsets.
Since each vertex in $S$ has at least $\lfloor t \rfloor$ neighbors,
each of these vertices can potentially connect to any of the $\ell$ subsets.

Given that $|S| > s \cdot \ell$, this would mean that on average,
each $\lfloor t \rfloor$-vertex subset of $A_1 \cup A_2$
must connect to more than $s$ vertices in $S$.
By the pigeonhole principle,
there must exist at least one specific $\lfloor t \rfloor$-vertex subset of
$A_1 \cup A_2$ that connects to more than $s$ vertices in $S$,
implying that this subset has over $s$ common neighbors in $S$.
Note that $u\notin A_1 \cup A_2$ and $S\subseteq B_1\subseteq N(u)$,
then $G$ contains a bipartite subgraph isomorphic to $K_{s,\lfloor t\rfloor+1}$.
Since $G$ satisfies the property $P_{t,r}$, we have the inequality
\[e(K_{s,\lfloor t\rfloor+1})=s(\lfloor t\rfloor+1)\leq t(s+\lfloor t\rfloor+1)+r.\]
Solving for $s$ in the inequality above yields that
$s\leq \frac{t(\lfloor t\rfloor+1)+r}{\lfloor t\rfloor+1-t },$
which contradicts the definition of $s$.
Therefore, it must be that $|S|\leq s\cdot \ell$.

Given that $p \leq |V_1| < (5s)^6$, it follows that $\ell$ and $|S|$ are constants.
For sufficiently large $n$, we can establish the following bounds.

First, consider the number of edges between $S$ and $A_1 \cup A_2$. We have
\[ e(S, A_1 \cup A_2) \leq |S|(|A_1| + |A_2|) \leq |S||V_1| \leq (5s)^{-6}n. \]
On the other hand, by the definition of $S$,
any vertex in $B_1 \setminus S$ has fewer than $\lfloor t \rfloor$ neighbors in $A_1 \cup A_2$.
Therefore,
\[ e(B_1 \setminus S, A_1 \cup A_2) \leq
|B_1 \setminus S|(\lfloor t \rfloor-1)\leq|N(u)|(\lfloor t\rfloor-1). \]
Combining these two inequalities, we obtain that
\[e(B_1, A_1 \cup A_2) \leq |N(u)|(\lfloor t \rfloor - 1) + (5s)^{-6}n.\]
It follows that
\begin{equation}\label{e5}
\sum_{v \in B_1, w \in N_{A_1 \cup A_2}(v)} x_w \leq
e\big(B_1, A_1 \cup A_2\big) x_{u^*}\leq
\Big(|N(u)|(\lfloor t \rfloor - 1) + (5s)^{-6}n\Big)x_{u^*}.
\end{equation}
Note that $x_u\leq x_{u^*}$.
Substituting inequality (\ref{e5}) into (\ref{a1}) gives that
\begin{align*}
\lfloor t \rfloor(n - \lfloor t \rfloor)x_u &\leq \Big(|N(u)|+\frac{(2s+1)n}{(5s)^6} +
\frac{2sn}{(5s)^4} + |N(u)|(\lfloor t \rfloor - 1)\Big)x_{u^*}\\
&\leq \lfloor t \rfloor \Big(|N(u)| + \frac{3sn}{(5s)^4} \Big)x_{u^*}.
\end{align*}
Rearranging terms yields that
$|N(u)|\geq(n-\lfloor t\rfloor)\frac{x_u}{x_{u^*}}-\frac{3sn}{(5s)^4}.$
For sufficiently large $n$, since $\lfloor t\rfloor\frac{x_u}{x_{u^*}}
\leq\lfloor t\rfloor<\frac{2sn}{(5s)^4}$, we conclude that
$|N(u)|\geq(\frac{x_u}{x_{u^*}}-\frac{1}{(5s)^3})n,$ as desired.
\end{proof}

In the following, we give the proof of Theorem \ref{t1}.
Recall that $V=\{v\in V(G):\ x_v\geq (5s)^{-1} x_{u^*}\}$. By the definition of $V$,
it is clear that $u^*\in V\subseteq V_1$.
We will show that $V$ is the desired set in Theorem \ref{t1}.
\vspace{1mm}

{\bf Proof of Theorem \ref{t1}:}
We start by proving lower bounds on $x_u$ and $|N(u)|$ for each vertex $u\in V$.
Suppose, for the sake of contradiction, that there exists $u_0 \in V$ such that
$x_{u_0}<(1-\frac{1}{50s^2})x_{u^*}.$

Since $u_0\in V$, we have $ x_{u_0} \geq \frac{x_{u^*}}{5s}.$
By Claim \ref{c6}, we know that
\[ |N(u^*)| \geq \Big(1-\frac{1}{(5s)^3}\Big)n \quad \text{and} \quad |N(u_0)|
\geq \Big(\frac{1}{5s}-\frac{1}{(5s)^3}\Big)n. \]

Define the following four sets: $A_1^* = V_1 \cap N(u^*)$, $B_1^* = N(u^*) \setminus A_1^*$,
 $A_2^* = V_1 \cap N_2(u^*)$ and $B_2^* = N_2(u^*) \setminus A_2^*$.
From Claim \ref{c5}, we have
$ |V_1| < (5s)^6 < \frac{n}{(5s)^3}$
for sufficiently large $n$. Thus,
\[ |B_1^*| \geq |N(u^*)| - |V_1| \geq \Big(1 - \frac{2}{(5s)^3}\Big)n. \]
Therefore,
\begin{align}\label{e6}
\big|B_1^* \cap N(u_0)\big|  \geq \big|B_1^*\big| + \big|N(u_0)\big| - n
 \geq \Big(\frac{1}{5s} - \frac{3}{(5s)^3}\Big)n
 \geq \frac{22n}{125s}.
\end{align}
Given (\ref{e6}), $u_0$ must have neighbors in $B_1^*$. Since $B_1^* \subseteq N(u^*)$,
$u_0$ is within distance two from $u^*$, implying $u_0 \in N(u^*) \cup N_2(u^*)$.
Recalling that $u_0 \in V \subseteq V_1$, we conclude $ u_0 \in A_1^* \cup A_2^*$.
Using inequality (\ref{a1}) with $u = u^*$, we get
\begin{align*}
\lfloor t\rfloor(n\!-\!\lfloor t \rfloor)x_{u^*} &\!\leq\!\Big(|N(u^*)|\!+\!\frac{2sn}{(5s)^6}
\!+\!\frac{2sn}{(5s)^4}\!+\!e\big(B_1^*,(A_1^*\cup A_2^*)\setminus\{u_0\}\big)\Big)x_{u^*}\!+\!e(B_1^*,\{u_0\})x_{u_0} \\
&\!\leq\!\Big(|N(u^*)|\!+\!\frac{2sn}{(5s)^6}\!+\!\frac{2sn}{(5s)^4}\!+\!e\big(B_1^*,A_1^* \cup A_2^*\big)\Big)x_{u^*}
\!+\!e(B_1^*,\{u_0\})(x_{u_0}\!-\!x_{u^*}),
\end{align*}
where $x_{u_0} - x_{u^*} < -\frac{x_{u^*}}{2(5s)^2}$ by our assumption.

From (\ref{e5}), we know that
$ e(B_1^*, A_1^* \cup A_2^*) \leq (\lfloor t \rfloor - 1)|N(u^*)| + \frac{n}{(5s)^6}.$
Furthermore, for large enough $n$, it easy to see
$\lfloor t \rfloor^2 \leq \frac{n}{(5s)^6}. $
This leads to
\begin{align*}
\lfloor t \rfloor n &\leq \lfloor t \rfloor |N(u^*)| + \frac{2sn}{(5s)^4}
+ \frac{(2s+2)n}{(5s)^6} - \frac{e(B_1^*, \{u_0\})}{2(5s)^2} \\
&< \lfloor t \rfloor n+\frac{2.2sn}{(5s)^4}-\frac{e(B_1^*, \{u_0\})}{2(5s)^2},
\end{align*}
which implies that
$e(B_1^*, \{u_0\}) < \frac{22n}{125s},$
contradicting (\ref{e6}).
Therefore, $x_u \geq(1-\frac{1}{50s^2})x_{u^*}$ for each $u\in V$.
By Claim \ref{c6}, this also means
\[|N(u)|\geq\Big(1-\frac{1}{2(5s)^2}-\frac{1}{(5s)^3}\Big)n\geq\big(1-\frac{1}{(5s)^2}\big)n.\]

To prove that $|V| = \lfloor t \rfloor$,
we begin by assuming that $|V| \geq \lfloor t \rfloor + 1$.
Based on the previous inequality,
every vertex $u \in V$ has at most $\frac{n}{(5s)^2}$ non-neighbors.
Hence, any $\lfloor t \rfloor + 1$ vertices in $V$ have at least $n-(\lfloor t \rfloor + 1)\frac{n}{(5s)^2}$
common neighbors, where $n-(\lfloor t \rfloor + 1)\frac{n}{(5s)^2}\geq s$ as $s>t+1$.
Then $G$ contains a subgraph isomorphic to $K_{\lfloor t \rfloor + 1, s}$.
By a similar argument to that in the proof of Claim \ref{c6},
there is a contradiction since $K_{\lfloor t\rfloor+1,s}$ does not satisfy the property $P_{t,r}$.
Therefore, $|V| \leq \lfloor t \rfloor.$

Next, suppose for contradiction that $|V| \leq \lfloor t \rfloor-1.$
Since $u^* \in V \setminus (A_1^* \cup A_2^*)$, we have
$ |V \cap (A_1^* \cup A_2^*)| \leq \lfloor t \rfloor - 2.$
This implies that
\[ e\big(B_1^*, V \cap (A_1^* \cup A_2^*)\big)\leq(\lfloor t\rfloor-2)n.\]
On the other hand, considering the edges between $B_1^*$ and $(A_1^*\cup A_2^*)\setminus V$,
we have
\[ e\big(B_1^*, (A_1^* \cup A_2^*) \setminus V\big) \leq e(G)\leq tn+r<sn.\]
By the definition of $V$, we know that
$x_w < \frac{x_{u^*}}{5s}$ for each $w \in (A_1^* \cup A_2^*) \setminus V$.
Now, setting $u = u^*$ in inequality (\ref{a1}), we obtain
\begin{align*}
\lfloor t\rfloor(n\!-\!\lfloor t\rfloor)x_{u^*}
&\leq \Big(|N(u^*)|
\!+\!\frac{2.5sn}{(5s)^4}\!+\!e\big(B_1^*,(A_1^*\cup A_2^*)\cap V\big)\Big)x_{u^*}
\!+\!e\big(B_1^*,(A_1^*\cup A_2^*)\setminus V\big)\frac{x_{u^*}}{5s}\\
&\leq \Big(n\!+\!\frac{n}{5}
\!+\!(\lfloor t\rfloor\!-\!2)n\!+\!\frac{n}{5}\Big)x_{u^*} \\
&= \big(\lfloor t\rfloor\!-\!\frac{3}{5}\big)nx_{u^*}.
\end{align*}
Dividing both sides by $x_{u^*}$ and rearranging terms, we find that
\[\lfloor t\rfloor(n-\lfloor t\rfloor)\leq\big(\lfloor t\rfloor-\frac{3}{5}\big)n.\]
Simplifying further, we have
$\lfloor t \rfloor^2 \geq \frac{3}{5}n$, a contradiction.
Therefore,
$ |V| = \lfloor t \rfloor$.

Combining both parts, we establish the desired results:
$x_u\ge\big(1-\frac{1}{50s^2}\big)x_{u^*}$
and $|N(u)|\ge \big(1-\frac{1}{25s^2}\big)n$ for each vertex $u\in V$.
Moreover, $|V|=\lfloor t\rfloor$.
Thus, we complete the proof of Theorem \ref{t1}. \qed

\section{Proof of Theorem \ref{t2}}

Let $G$ be a graph that achieves the maximal spectral radius in $\mathcal{G}$.
It is not hard to verify that $K_{\lfloor t\rfloor,n-\lfloor t\rfloor}\in \mathcal{G}$.
Therefore,
$\rho(G)\geq\rho(K_{\lfloor t\rfloor,n-\lfloor t\rfloor})
=\sqrt{\lfloor t\rfloor (n-\lfloor t\rfloor)}$.

Let $V \subseteq V(G)$ be a set defined in accordance with Theorem \ref{t1}.
We define $U$ as the set of vertices that are complete to $V$,
meaning every vertex in $U$ is adjacent to every vertex in $V$.
It is clear that $U$ and $V$ are disjoint sets, i.e.,
$U \cap V = \emptyset$. We then define $W = V(G) \setminus (U \cup V)$.

\begin{claim}\label{c1}
 $|U|\geq \frac{n}{2}$.
\end{claim}
\begin{proof}
Set $V=\{v_1,v_2,\dots,v_{\lfloor t\rfloor}\}$.
By the Inclusion-Exclusion Principle, we have
\begin{equation*}
|U|=\big|\cap_{i=1}^{\lfloor t\rfloor}N(v_i)\big|
\geq\sum\limits_{i=1}^{\lfloor t\rfloor}\big|N(v_i)\big|-
(\lfloor t\rfloor-1)\big|\cup_{i=1}^{\lfloor t\rfloor}N(v_i)\big|.
\end{equation*}
It follows that
\[ |U|\geq\Big(1-\frac{1}{25s^2}\Big)n\cdot\lfloor t\rfloor-(\lfloor t\rfloor-1)n
=\Big(1-\frac{\lfloor t\rfloor}{25s^2}\Big)n
\geq\frac{n}{2},\]
as claimed.
\end{proof}

\begin{claim}\label{c2}
 $W=\emptyset$.
\end{claim}

\begin{proof}
Choose a vertex $u\in U\cup W$.
By Claim \ref{c1}, let $N_{U}(u)=\{u_1,u_2,\dots,u_d\}$ and
denote by $H$ the subgraph in $G$ induced by the set $V\cup N_U(u)\cup \{u\}$.
Then $e(H)\geq (\lfloor t\rfloor+1)d$. Since $H $ has the property $P_{t,r}$,
we have $e(H)\leq t|H|+r=t(\lfloor t\rfloor+1+d)+r$.
Combining the two inequalities, we have
$(\lfloor t\rfloor+1-t)d\leq t(\lfloor t\rfloor+1)+r$, and so
\begin{equation}\label{du}
d_U(u)=d\leq \frac{t(\lfloor t\rfloor+1)+r}{\lfloor t\rfloor+1-t }\leq s-1.
\end{equation}
for each vertex $u\in U\cup W.$

Note that $W$ is not complete to $V$, meaning each vertex in $W$ is not adjacent to some vertex in $V$.
By Theorem \ref{t1}, every vertex in $V$ is not adjacent to at most $\frac{n}{25s^2}$ vertices of $G$.
Recall that $s>t+1$. Therefore,
\begin{equation}\label{duu}
|W|\leq \frac{\lfloor t\rfloor n}{25s^2}\leq \frac{n}{25s}.
\end{equation}
Note that $d(u)\leq d_U(u)+|V|+d_W(u)\leq s-1+\lfloor t\rfloor+d_W(u)$.
Moreover, since $G$ has the property $P_{t,r}$,
either $e(W)\leq t|W|+r$ or $|W|<t$.
Consequently, we get that
\[\sum\limits_{w\in W}\rho x_w
\leq \Big(\big(s-1+\lfloor t\rfloor\big)|W|+2e(W)\Big)x_{u^*}
\leq \frac{(s+3t)\lfloor t\rfloor n}{25s^2}x_{u^*}.\]
It follows that
\begin{equation}\label{e2}
\sum\limits_{w\in N_W(v)} x_w\leq\sum\limits_{w\in W} x_w
< \frac{(s+3t)\lfloor t\rfloor n}{25\rho s^2}x_{u^*}.
\end{equation}
Now for $u\in U\cup W$, we have $\rho x_u\leq(s-1+\lfloor t\rfloor)x_{u^*}+\sum_{w\in N_W(u)} x_w $.
Combining with the inequality \eqref{e2} and
$\rho\geq\sqrt{\lfloor t\rfloor (n-\lfloor t\rfloor)}$,
we deduce that
\begin{equation}\label{e3}
x_u<\Big(\frac{s-1+\lfloor t\rfloor}{\rho}+\frac{(s+3t)\lfloor t\rfloor n}{25\rho^2s^2}\Big)x_{u^*}
\leq \frac{s+4t}{25s^2}x_{u^*}<\frac{x_{u^*}}{5s}.
\end{equation}

Suppose that $W\neq\emptyset$.
Let $G^*$ denote the graph obtained from $G$ by
removing all edges within $W$ and from $W$ to $U$,
then adding all possible edges from $W$ to $V$.
Additionally, let $w_1$ be a vertex of minimal degree in $G[W]$,
and let $E_1$ be the set of edges incident to  $w_1$ in $G[W]$.
Following the same approach, let $w_2$ be a vertex of minimal degree
in $G[W\setminus \{w_1\}]$,
and let $E_2$ be the set of edges incident to $w_2$ in $G[W\setminus \{w_1\}]$.
Repeating the process,
we define $W_1=W$ and $W_i=W\setminus\{w_1,\dots,w_{i-1}\}$ for $i\geq2$,
where $w_i$ is a vertex of minimum degree in $G[W_i]$
and $E_i$ is the set of edges that are incident to $w_i$ in $G[W_i]$.
By continuing this process iteratively, we have
\[E(G[W])=\cup_{i=1}^{|W|}E_i=\cup_{i=1}^{|W|}\{w_iw:~ w\in N_{W_i}(w)\}.\]
Recall that $G\in \mathcal{G}$ and
$s=\big\lceil \frac{t(\lfloor t\rfloor+1)+r}{\lfloor t\rfloor+1-t}\big\rceil+1\geq t^2+r+1$.
When $|W_i|\geq t$, we have
$$d_{W_i}(w_i)\leq\frac{2e(W_i)}{|W_i|} \leq\frac{2t|W_i|+2r}{|W_i|}\leq 2t+\frac{2r}{t}<2s.$$
When $|W_i|<t$, we also get
$d_{W_i}(w_i)\leq\frac{2\binom{|W_i|}{2}}{|W_i|}<t<2s.$
Now, we can drive that
\begin{align}\label{e4}
\rho(G^*)-\rho(G)&\geq X^T\big(A(G^*)-A(G)\big)X \nonumber\\
&=\sum\limits_{i=1}^{|W|}2x_{w_i}\Big(\sum\limits_{v\in V\setminus N_V(w_i)}x_v
-\sum\limits_{u\in N_U(w_i)}x_u-\sum\limits_{w\in N_{W_i}(w_i)}x_w\Big).
\end{align}
In view of the definition of $W$,
we have $|V\setminus N_V(w_i)|\geq1$ for each $w_i\in W$.
By Theorem \ref{t1}, $x_v\ge\big(1-\frac{1}{50s^2}\big)x_{u^*}$ for each $v\in V$.
Moreover, by inequality \eqref{du}, we know that
$d_U(u)\leq s-1$ for each $u\in U$, and based on the previous discussion,
we have $d_{W_i}(w_i)<2s$ for each $w_i\in W_i$.
Combining \eqref{e3} and \eqref{e4}, we can see that
\[\rho(G^*)-\rho(G)\geq \sum\limits_{w_i\in W}2x_{w_i}\Big(\big(1-\frac{1}{50s^2}\big)
-\big(3s-1\big)\frac{1}{5s}\Big)x_{u^*}>0.\]

Based on the maximality of $\rho(G)$, the new graph $G^*$ cannot satisfy the property $P_{t,r}$.
Specifically, there exists an induced subgraph $H$ of $G^*$ such that
$|H| \geq t$ and $e(H) > t|H| + r$. Since $G^*[U \cup V] = G[U \cup V]$,
the violation of the property $P_{t,r}$ must involve vertices in $W$,
i.e., $V(H) \cap W \neq \emptyset$.

Let $H' = G^*[V(H) \setminus W]$. According to the definition of $G^*$,
we know that $N_{G^*}(w)=V$ for each $w\in W$. Hence,
$e(H)=e(H')+|V(H)\cap W|\cdot|V(H)\cap V|.$
Since $|V(H) \cap V| \leq \lfloor t \rfloor \leq t$ and $e(H) > t|H| + r$,
it follows that $$e(H') > t|V(H)| + r - t|V(H) \cap W| = t|V(H')| + r.$$

If $|V(H')|\geq t$, the above inequality leads to a contradiction
because $G^*[V(H')] = G[V(H')]$
and $G$ is assumed to satisfy the property $P_{t,r}$.

If $|V(H')|< t$, then by Claim \ref{c1} and the inequality in (\ref{duu}),
we have $|U\setminus V(H')|\geq \frac n2-t\geq |V(H)\cap W|$.
Since $U\setminus V(H)=U\setminus V(H')$,
we can choose a subset $U'$ of $U\setminus V(H)$ with $|U'|=|V(H)\cap W|$.
Denote $H'' = G^*[U'\cup(V(H)\setminus W)]$.
Then $|V(H'')|=|V(H)|$, and
$N_{G^*}(w)=V\subseteq N_{G^*}(u)$ for every $w\in W$ and every $u\in U'$.
Thus, $e(H'')\geq e(H)>t|V(H)|+r=t|V(H'')|+r$.
This also leads to a contradiction.

In both cases, we arrive at a contradiction,
which implies that our initial assumption $W \neq \emptyset $ must be false.
Therefore, $ W = \emptyset $.
\end{proof}

By Claim \ref{c2}, we have $V(G) = V \cup U$ and $e(V,U)=|V||U|$,
where $|V|=\lfloor t \rfloor$. Based on \eqref{du},
every $u \in U$ has at most $s-1$ neighbors in $U$.
Let $u^\star\in U $ such that $x_{u^\star}=\max_{u\in U}x_u$.
Then
$\rho x_{u^\star}\leq (s-1)x_{u^\star} + \lfloor t \rfloor x_{u^*}.$
It follows that
$x_{u^\star} \leq \frac{\lfloor t \rfloor}{\rho - s+1} x_{u^*}.$
Also, since $u^* \in V$, we have
$\rho x_{u^*} \leq (\lfloor t \rfloor-1)x_{u^*} + (n - \lfloor t \rfloor)x_{u^\star}.$
Combining these two inequalities, we get that
\[(\rho-\lfloor t \rfloor+1) x_{u^*} \leq (n - \lfloor t \rfloor) \frac{\lfloor t \rfloor}{\rho-s+1} x_{u^*}. \]
Dividing both sides by $x_{u^*}$, we obtain
$\rho-\lfloor t \rfloor+1\leq\frac{(n - \lfloor t \rfloor) \lfloor t \rfloor}{\rho-s+1}.$
Expanding and rearranging terms, we have
\[ \rho^2 - (s + \lfloor t \rfloor-2) \rho - \lfloor t \rfloor n
+\lfloor t \rfloor^2+(s-1)(\lfloor t \rfloor-1) \leq 0. \]
Solving this quadratic inequality, we find that
\[ \rho \leq \frac{(s + \lfloor t \rfloor-2)}{2} + \sqrt{\frac{(s + \lfloor t \rfloor-2)^2}{4}
+ \lfloor t \rfloor n-\lfloor t \rfloor^2-(s-1)(\lfloor t \rfloor-1)}. \]
Therefore,
$\rho \leq c(s, t) + \sqrt{\lfloor t \rfloor n},$
where $c(s, t)$ is a constant depending on $s$ and $t$. This completes the proof
of Theorem \ref{t2}.

\section{Characterization of extremal graphs}
In this section, we will present the proof Theorem \ref{t3}.
Before delving into our proof, we introduce some concepts and several lemmas.
An $\ell$-walk consists of $\ell + 1$ vertices where consecutive vertices are adjacent.
Let $w_G^\ell(u)$ be the number of $\ell$-walks starting at vertex $u$ in $G$.
The total number of $\ell$-walks in $G$
is defined as $W_\ell(G) = \sum_{v \in V(G)} w_G^\ell(v)$.
Obviously, $W_1(G)=2e(G)$ and $W_2(G) = \sum_{v \in V(G)}(d_G(v))^2.$

Let $G_1,G_2$ be two graphs, potentially of different orders.
We say that $G_1$ is {\it walk-equivalent} to $G_2$ and denote this as $G_1\equiv G_2$,
if $W_\ell(G_1) = W_\ell(G_2)$ for all $\ell\geq 1$.
If there exists an $\ell$ such that $W_\ell(G_1) > W_\ell(G_2)$
and $W_i(G_1) = W_i(G_2)$ for $i \leq \ell - 1$,
we say that $G_1$ is {\it strictly walk-preferable} to $G_2$ and denote this as $G_1 \succ G_2$.

The following theorem, provided by Zhang,
offers a framework for comparing the spectral radius of graphs based on their local walk structures.
This seems particularly useful in various graph-theoretic analyses and applications.

\begin{lemma}(\cite{zhang24})\label{l0}
Let $G_0$ be an $n$-vertex connected graph. 
Assume that $S$ is a non-empty proper subset of $V(G_0)$
and $T$ is a set of isolated vertices in $G_0-S$ such that $e(S,T)=|S||T|$.
For $i\in\{1,2\}$, let $H_i$ be a graph with vertex set $T$,
and let $G_i$ be the graph obtained from $G_0$ by embedding the edges of $H_i$ into $T$.
When $\rho(G_0)$ is sufficiently large compared to $|T|$, we have the following conclusions:\vspace{1mm}\\
(i) If $H_1 \equiv H_2$, then $\rho(G_1) = \rho(G_2)$.\\
(ii) If $H_1 \succ H_2$, then $\rho(G_1) > \rho(G_2)$.
\end{lemma}

Let $\mathbb{G}(n, m)$ be the family of graphs with $n$ vertices and $ m $ edges,
and let $ F_p(H) $ be the sum of the $ p $-th powers of vertex degrees of a graph $ H $,
formally defined as
\[
F_p(H) = \sum_{v \in V(H)}\big(d_H(v)\big)^p,
\]
where $ p \geq 1 $ is an integer.
Ismailescu and Stefanica showed the following result regarding the maximization of $F_p(H)$.

\begin{lemma}\label{l1}(\cite{Ismailescu})
 Let $H^* \in \mathbb{G}(n, m)$ such that $F_p(H^*) = \max_{H \in\mathbb{G}(n, m)} F_p(H)$,
 where $p \geq 2 $ is an integer. If $ n \geq m + 2 $,
 then $ H^* $ is the star $ K_{1,m} $ plus $ n - m - 1 $ isolated vertices.
 Moreover, it is unique except in the case when $ p = 2 $ and $ m = 3 $
 (when both $ K_{1,3} \cup (n - 4)K_1 $ and $ K_3 \cup (n - 3)K_1 $ are extremal graphs).
\end{lemma}

Given that $p=2$ and $n=m-1$,
and following a similar approach to Lemma \ref{l1},
the result below can be derived.

\begin{lemma}\label{l2}
Let $H^* \in \mathbb{G}(n, m)$ such that $F_2(H^*) = \max_{H \in\mathbb{G}(n, m)} F_2(H)$.
If $m = n - 1$, then $H^*\cong K_{1,m}$. Moreover, this configuration is unique unless $m = 3$,
in which case both $K_{1,3}$ and $K_3 \cup K_1$ are extremal graphs.
\end{lemma}

Given a graph $F$,
we define $\eta(F)=e(F)+(\lfloor t\rfloor-t)|V(F)|$, 
and let $\mathcal{C}(F)$ be the set of all components of $F$.
It is easy to see that $\eta(F)=\sum_{F_i\in \mathcal{C}(F)}\eta(F_i)$.
Furthermore, a component $F_i$ of $F$ is said to be {\it positive} if $\eta(F_i)>0$.
Let $\mathcal{C}^+(F)$ be the set of all positive components of $F$.
We denote
\[\eta^+(F)=\sum_{F_i\in \mathcal{C}^+(F)}\eta(F_i),
\ t_0=\lfloor t\rfloor+1-t\ \text{and}\ g_0=t\lfloor t\rfloor+r-\binom{\lfloor t\rfloor}{2}.\]

Specially, if $\mathcal{C}^+(F)=\emptyset$, then $\eta^+(F)=0$. 
Recall that $r\geq-\binom{\lfloor t\rfloor+1}{2}$.
Clearly, $t_0$ and $g_0$ are two constants with $0<t_0\leq1$ and $g_0\geq0$.
Based on the definitions of $g_0$ and $\eta(F)$, we have the following two lemmas.

\begin{lemma}\label{f-1}
If $F'$ is a subgraph of $F$, then $\eta(F')\leq\eta(F)$.
\end{lemma}

\begin{proof}
We first assume that $F$ is connected. 
If $V(F')=V(F)$, then it is clear that $\eta(F')\leq\eta(F)$.
Next, suppose that $V(F')\subsetneq V(F)$,
and let $F_1',\dots,F_p'$ be the components obtained by removing $V(F')$ from $F$. 
Since each $F_i'$ is connected, $e(F_i') \geq |V(F_i')| - 1$, 
and there is at least one edge connecting $F_i'$ to $F'$. Thus
 \begin{align*}
   \eta(F) &= e(F) + (\lfloor t \rfloor - t)|V(F)| \\
   &\geq e(F') + \sum_{i=1}^p e(F_i') + \sum_{i=1}^p e(F_i',F')  
   + (\lfloor t \rfloor - t)|V(F')| + (\lfloor t \rfloor - t)\sum_{i=1}^p |V(F_i')|.
 \end{align*}
Combining $e(F_i') \geq |V(F_i')| - 1$ and $e(F_i',F')\geq 1$, we have
\[
\eta(F) \geq \eta(F') + \sum_{i=1}^p (|V(F_i')| - 1) + k 
+ (\lfloor t \rfloor - t)\sum_{i=1}^p |V(F_i')| = \eta(F') + t_0 \sum_{i=1}^p |V(F_i')|,
\]
where $t_0 =\lfloor t \rfloor+ 1 - t> 0$. It follows that $\eta(F')\leq\eta(F)$.

Now we consider the case when $F$ is not connected.
Let $\mathcal{C}(F)=\{F_1,\dots,F_q\}$. 
Based on the above argument, 
we know that $\eta\big(F'[V(F')\cap V(F_i)]\big)\leq\eta(F_i)$ for $i\in\{1,\ldots,q\}$.
It follows that $\eta(F')\leq\sum_{i=1}^q\eta(F_i)=\eta(F)$, as desired.
\end{proof}

\begin{lemma}\label{f1}
$K_{\lfloor t\rfloor}\nabla F$ satisfies the property $P_{t,r}$ if and only if
$\eta^+(F)\leq g_0$.
\end{lemma}

\begin{proof}
First, assume that $K_{\lfloor t \rfloor} \nabla F$ satisfies the property $P_{t,r}$. 
Recall that $g_0=t\lfloor t\rfloor+r-\binom{\lfloor t\rfloor}{2}\geq0$.
If $\mathcal{C}^+(F)=\emptyset$, then $\eta^+(F)=0\leq g_0$, 
and the conclusion holds trivially. Otherwise, let $F^+$ be the union of all positive components of $F$. 
Since $K_{\lfloor t \rfloor} \nabla F^+$ is an induced subgraph of $K_{\lfloor t \rfloor} \nabla F$ 
with $|V(K_{\lfloor t \rfloor} \nabla F^+)| \geq t$, it must satisfy $P_{t,r}$, and hence,
\[
e(K_{\lfloor t \rfloor} \nabla F^+) = e(F^+) + \lfloor t \rfloor |V(F^+)| 
+ \binom{\lfloor t \rfloor}{2} \leq t\big(\lfloor t \rfloor + |V(F^+)|\big) + r.
\]
Rearranging this inequality, we obtain that
$e(F^+)+(\lfloor t\rfloor-t)|V(F^+)|\leq g_0,$
which implies that $\eta^+(F)=\eta(F^+)\leq g_0$.

Now, assume that $\eta^+(F) \leq g_0$. 
We will prove that $K_{\lfloor t\rfloor}\nabla F$ satisfies the property $P_{t,r}$.
Let $H$ be an arbitrary subgraph of $K_{\lfloor t \rfloor} \nabla F$ with $|V(H)| \geq t$. 
It suffices to show that $e(H)\leq t|V(H)|+r.$
Denote $F'=F[V(H)\cap V(F)]$ and $t'=|V(H)\cap V(K_{\lfloor t \rfloor})|$,
where $0\leq t'\leq\lfloor t \rfloor$.
By Lemma \ref{f-1}, we have
$\eta(F')\leq\eta(F) \leq \eta^+(F) \leq g_0,$
and hence
\begin{align}\label{ee4}
e(F')=\eta(F')-(\lfloor t\rfloor-t)|V(F')|\leq g_0-(\lfloor t\rfloor-t)(|V(H)|-t').
\end{align}
Observe that $e(H)\leq e(F')+\binom{t'}{2}+t'(|V(H)|-t')$, and recall that
$g_0=t\lfloor t\rfloor+r-\binom{\lfloor t\rfloor}{2}$.
By combining (\ref{ee4}), we obtain that
   \begin{align*}
      e(H)\leq t|V(H)| + r-(\lfloor t\rfloor-t')(|V(H)|-t'-t)
     +\binom{t'}{2}- \binom{\lfloor t\rfloor}{2}.
   \end{align*}
Note that $\lfloor t\rfloor \geq t'$ and $|V(H)|\geq t$. 
Define $f(t')=(\lfloor t\rfloor -t')t'+\binom{t'}{2}-\binom{\lfloor t\rfloor}{2}$. 
Then
\begin{align*}
e(H)\leq t|V(H)|+r+f(t').
\end{align*}
If $t' = \lfloor t \rfloor$, then $f(t')=0$ and thus $e(H) \leq t|V(H)| + r$, as required. 
If $t'\leq \lfloor t \rfloor - 1$, it is easy to check that $f(\lfloor t\rfloor-1)=0$ and $f'(t')=\lfloor t\rfloor-t'-\frac12>0$.
Thus $f(t') \leq 0$ for $0\leq t' \leq \lfloor t \rfloor - 1$.
Hence, $e(H) \leq t|V(H)| + r$.
This completes the proof.
\end{proof}

Let $G$ be a graph that attains the maximum spectral radius over all
$n$-vertex graphs satisfying the property $P_{t,r}$.
By Claim \ref{c2}, we know that $V(G)=V\cup U$, where $|V|=\lfloor t\rfloor$
and $e(V,U)=|V||U|$.
Let $T$ be a subset of $U$
such that $e(T,U\setminus T)=0$ and $|T|$ is constant.
Define two new graphs $H_1$ and $H_2$
with $V(H_1) = V(H_2) = T$ and $e(H_1) = e(H_2)$.
Let $G_i$ be the graph obtained from $G$ by embedding the edges of $H_i$ into $T$ for $i \in \{1, 2\}$.
By Lemma \ref{l0}, we obtain the following corollary.

\begin{cor}\label{f0}
If $H_1\succ H_2$ (specially, $F_2(H_1)>F_2(H_2)$), then $\rho(G_1)>\rho(G_2)$.
\end{cor}

Let $g(V)=t\lfloor t\rfloor+r-e(V)$. 
It is obvious that $g_0\leq g(V)\leq t\lfloor t\rfloor+r$. 
Observe that $G[V]\nabla F$ is a spanning subgraph of $K_{\lfloor t\rfloor}\nabla F$.
By a similar argument as in the first part of Lemma \ref{f1}, we obtain the following statement.

\begin{cor}\label{f2}
If $ G[V]\nabla F$ satisfies the property $P_{t,r}$, then
$\eta^+(F)\leq g(V).$ 
\end{cor}

We will now formally present the proof of Theorem \ref{t3},
starting with the following key claim.
The proof proceeds via a minimal counterexample.

\begin{claim}\label{cc1}
  $V$ is a clique.
\end{claim}

\begin{proof}
Suppose to the contrary that $G[V]$ has a non-edge $v_1v_2$.
Clearly, $\rho(G+v_1v_2)>\rho(G)$. Then
$G$ admits a subgraph $H$ of order at least $t$ such that $H+v_1v_2$ no longer satisfies the property $P_{t,r}$.
Since $H$ satisfies the property $P_{t,r}$ but $H+v_1v_2$ no longer does, 
we conclude that $e(H)=\big\lfloor t|V(H)|+r\big\rfloor$. 
Now, $H+v_1v_2$ is called a counterexample of $G+v_1v_2$.
We may assume that $H+v_1v_2$ is minimal with respect to the order.

Let $A = V \cap V(H)$, $B = U \cap V(H)$ and $C= V \setminus V(H)$.
Define $L$ as the vertex set of all components in $G[U\setminus B]$
that have at least one edge connecting to some vertex in $B$.
Furthermore, let $N_0$ denote the vertex set of all tree-components in $G[U \setminus (B \cup L)]$
and $N_1$ be the vertex set of all other components in $G[U \setminus (B \cup L)]$.

First, we will prove that the number of vertices in $L \cup N_1$ is finite.
Consider the size of the induced subgraph $G[V \cup B \cup L \cup N_1]$:
\begin{align*}
e(V \cup B \cup L \cup N_1) &=e(V\cup B)+e(L\cup N_1)+e(V\cup B,L\cup N_1)\\
& \geq e(H)+e(L)+e(N_1)+e(B,L)+e(V,L\cup N_1).
\end{align*}
Recall that $e(H)=\big\lfloor t|V(H)|+r\big\rfloor$. Moreover,
we see that $e(L)+e(B,L)\geq |L|$, $e(N_1) \geq |N_1|$ 
and $e(V, L \cup N_1)=\lfloor t\rfloor(|L|+|N_1|)$.
Thus
\[
e(V \cup B \cup L \cup N_1) \geq\big\lfloor t|V(H)|+r\big\rfloor
+ |L| + |N_1|+\big\lfloor t \big\rfloor\big(|L| + |N_1|\big).
\]
On the other hand, by the property $P_{t,r}$, we have
\[
e(V \cup B \cup L \cup N_1) \leq t(|V|+|B| + |L| + |N_1|) + r.
\]
Note that $|V|+|B|=|V(H)|+|C|$. Combining the above two inequalities, we obtain
\begin{align*}
\big(\lfloor t \rfloor + 1 - t\big)\big(|L| + |N_1|\big) &\leq t\big(|V(H)|+|C|\big)+r-\big\lfloor t|V(H)|+r\big\rfloor \\
&\leq t|C|+1\leq t\lfloor t \rfloor +1,
\end{align*}
which shows that $|L| + |N_1|$ is finite.

Recall that $t_0=\lfloor t\rfloor+1-t$.
We now partition $N_0$ into two parts:
Let $N_0^1$ denote the vertex set of all components in $G[N_0]$ with cardinality at least $\lfloor1/{t_0}\rfloor+1$,
and let $N_0^2$ denote the vertex set of all remaining components in $G[N_0]$.

Next, we will show that $|N_0^1|$ is also finite. 
Let $F$ be a component of $G[N_0^1]$. 
Note that $F$ is a tree with $|V(F)|\geq\lfloor1/{t_0}\rfloor+1>1/{t_0}$, 
and recall that 
$\eta(F)=e(F)+(\lfloor t\rfloor-t)|V(F)|$.
It follows that $\eta(F)=t_0|V(F)|-1>0$, 
and thus, $\eta^+(N_0^1)=\eta(N_0^1)=t_0|N_0^1|-c$, 
where $c$ is the number of the components in $G[N_0^1]$.
Notice that $c(\lfloor1/{t_0}\rfloor+1)\leq|N_0^1|$.
Let $t'_0=1/(\lfloor1/{t_0}\rfloor+1)$. 
Then $ t'_0<t_0$, $c\leq t'_0|N_0^1|$ and
$\eta^+(N_0^1)\geq(t_0-t_0')|N_0^1|$.
Since $G[V]\nabla G[N_0^1]$ 
satisfies the property $P_{t,r}$, by Corollary \ref{f2}, 
we get that $\eta^+(N_0^1)\leq g(V)\leq t\lfloor t\rfloor+r$. 
By combining these two inequalities, the assertion follows.

Thirdly, we shall demonstrate that $|B|$ is finite provided that $C\neq\emptyset$.
Consider the subgraph $G[V(H)\cup C]$.
It is easy to see that
$e(V(H)\cup C)\geq e(H)+|C||B|$,
where $e(H)=\big\lfloor t|V(H)|+r\big\rfloor\geq t|V(H)|+r-1$.
However, using the property $P_{t,r}$, we have
$e(V(H)\cup C)\leq t(|V(H)|+|C|)+r.$
Hence, $|C||B| \leq t|C| + 1,$
and so $|B| \leq t + 1$.

Observe that $e(B,N_1\cup N_0^1\cup N_0^2)=0$.
In view of \eqref{du}, we have $e(B, L)\leq |L|(s-1)$.
Building on the previous argument, we further see that both
$e(B,L)$ and $e(L\cup N_1\cup N_0^1)$ are constants.
We proceed with the proof by considering two cases.

The first case is when $e(B)$ is finite.
Now, $e(U \setminus N_0^2)$ is also finite.
Let $G'$ be the graph obtained from $G$ by removing all edges within $U \setminus N_0^2$
and adding the edge $v_1v_2$.
We assert that $G'$ satisfies the property $P_{t,r}$. 
To verify this, observe that $E(G'[U]) = E(G'[N_0^2])$. 
Let $F'$ be a component of $G'[U]$. Then $F'$ is a tree with 
$|V(F')|\leq\lfloor\frac{1}{t_0} \rfloor$.
Thus 
\begin{align}\label{z00}
\eta(F')=e(F')+(\lfloor t\rfloor-t)|V(F')|=t_0|V(F')|-1\leq 0,
\end{align}
which implies that $G'[U]$ contains no any positive components.
Hence, $\eta^+(G'[U]) = 0 \leq g_0$. 
By Lemma \ref{f1}, $K_{\lfloor t \rfloor} \nabla G'[U]$ satisfies the property $P_{t,r}$. 
Since $G'$ is a spanning subgraph of $K_{\lfloor t \rfloor} \nabla G'[U]$, 
it follows that $G'$ also satisfies the property $P_{t,r}$. 

By Theorem \ref{t1}, $x_v \geq \big(1 - \frac{1}{50s^2}\big)x_{u^*}$ for $v \in V$.
In view of (\ref{du}), we have $x_u \leq \frac{\lfloor t \rfloor}{\rho -s+1} x_{u^*}$ for each $u \in U$.
Now we see that $$\rho(G')-\rho(G) \geq X^T\big(A(G')-A(G)\big)X=
2x_{v_1}x_{v_2}-\!\!\sum_{u_iu_j\in E(G[U\setminus N_0^2])}\!\!2x_{u_i}x_{u_j}.$$
Recall that $\rho\geq\sqrt{\lfloor t\rfloor(n-\lfloor t\rfloor)}$ and $e(G[U \setminus N_0^2])$ is finite.
Therefore, we can deduce that $\rho(G')>\rho(G)$, which contradicts the maximality of $\rho(G)$.

It remains the case when $e(B)$ is sufficiently large.
Now $B$ is an infinite set,
and by the previous proof, we have $C=\emptyset$. Thus $V(H)=V \cup B$.
We construct $G''$ as follows:
(i) delete all edges in $G[L \cup N_1 \cup N_0^1]$;
(ii) delete all edges from $L$ to $B$;
(iii) delete an edge in $G[B]$ and add the edge $v_1v_2$ to $G[V]$.
Using a similar argument to the one for $G'$,
we have $\rho(G'')>\rho(G)$,
which implies that $G''$ does not satisfy the property $P_{t,r}$.

We conclude the proof by contradiction, demonstrating that 
$G''$ satisfies the property $P_{t,r}$.
Observe that $e(G''[V \cup B])=e(G[V \cup B])$.
Since $G$ satisfies the property $P_{t,r}$,
we have $e(G''[V \cup B])=e(G[V \cup B])\leq t(|V|+|B|)+r$.
Recall that $V(H)=V \cup B$ and $H + v_1v_2$ is a minimal counterexample.
We see that $e(H')\leq t|V(H')|+r$ for any subgraph $H'$ of $G''[V \cup B]$.
Hence, $G''[V \cup B]$ satisfies the property $P_{t,r}$.

On the other hand,
observe that $G''[U\setminus B]\cong G'[U\setminus B]$,
and every component $F'$ of $G''[U\setminus B]$ is a tree with 
$|V(F')|\leq\lfloor\frac{1}{t_0} \rfloor$.
In view of (\ref{z00}),
we know that $\eta(F')\leq0$ for each component $F'$ of $G''[U\setminus B]$.
It follows that $\eta(G''[U\setminus B])\leq0$.

Let $H''$ be an arbitrary subgraph of $G''$.
Moreover, let $H_1''=G''[V(H'')\cap (V\cup B)]$ and $H_2''=G''[V(H'')\cap (U\setminus B)]$.
Since $G''[V \cup B]$ satisfies the property $P_{t,r}$,
we have $e(H_1'')\leq t|V(H_1'')|+r$.
By Lemma \ref{f-1}, $\eta(H_2'')\leq \eta(G''[U\setminus B])\leq0$,
which yields that $e(H_2'')+\lfloor t\rfloor|V(H_2'')|\leq t|V(H_2'')|.$
Note that there are no edges between $B$ and $U\setminus B$ in $G''$.
Thus, 
$e(H'')\leq e(H_1'')+e(H_2'')+\lfloor t\rfloor|V(H_2'')|\leq t|V(H'')|+r.$
Therefore, $G''$ satisfies the property $P_{t,r}$.
This completes the proof.
\end{proof}

By Lemma \ref{f1} and Claim \ref{cc1}, we have $\eta^+(G[U])\leq g_0$.
Next, we classify the components of $G[U]$ into three distinct categories
$\mathcal{N}_1,$ $\mathcal{N}_2$ and $\mathcal{N}_3$, where

(i) $\mathcal{N}_1$ includes all components that are not trees;

(ii) $\mathcal{N}_2$ contains all tree components of order at least $\lfloor1/t_0\rfloor+1$;

(iii) $\mathcal{N}_3$ consists of all tree components of order at most $\lfloor1/t_0\rfloor$.

\begin{claim}\label{cc2}
Every component in $\mathcal{N}_1$ is positive. 
Moreover, $\sum_{F\in \mathcal{N}_1} |V(F)| \leq g_0/t_0.$
\end{claim}

\begin{proof}
For each $F\in\mathcal{N}_1$,
we have $e(F) \geq |V(F)|$.
Recall that $t_0=\lfloor t\rfloor+1-t>0$.
Thus,
\begin{align}\label{z00-2}
\eta(F) = e(F) + (\lfloor t \rfloor- t)|V(F)|\geq t_0 |V(F)|>0.
\end{align}
Hence, $F$ is a positive component. 
Furthermore, summing over all components $F\in \mathcal{N}_1$, we obtain that
$\sum_{F\in \mathcal{N}_1} \eta(F)\geq t_0 \sum_{F\in \mathcal{N}_1} |V(F)|$.
However, by Lemma \ref{f1}, we know that $\sum_{F\in \mathcal{N}_1} \eta(F)\leq
\eta^+(G[U])\leq g_0$
since $G\cong K_{\lfloor t \rfloor}\nabla G[U]$ and $G$ satisfies the property $P_{t,r}$.
Therefore, we have
$\sum_{F\in \mathcal{N}_1} |V(F)| \leq g_0/t_0,$
as claimed.
\end{proof}

\begin{claim}\label{cc3}
Every component in $\mathcal{N}_2$ is positive and $\sum_{F\in \mathcal{N}_2} |V(F)|$ is constant. 
\end{claim}

\begin{proof}
For each $F\in\mathcal{N}_2$, we have $e(F)=|V(F)|-1$
and $|V(F)|\geq\lfloor1/t_0\rfloor+1$.
Thus, 
\begin{align}\label{z00-1}
\eta(F) = e(F) + (\lfloor t \rfloor- t)|V(F)|= t_0 |V(F)| - 1>0.
\end{align}
Hence, all components in $\mathcal{N}_2$ are positive and 
$\sum_{F\in \mathcal{N}_2} \eta(F)=t_0\sum_{F\in \mathcal{N}_2} |V(F)|-|\mathcal{N}_2|$,
where $\sum_{F\in \mathcal{N}_2} |V(F)|\geq |\mathcal{N}_2|(\lfloor1/t_0\rfloor+1)$.
Thus, $$\sum_{F\in \mathcal{N}_2} \eta(F)
\geq\Big(t_0-\frac1{\lfloor1/t_0\rfloor+1}\Big)\sum_{F\in \mathcal{N}_2} |V(F)|.$$
On the other hand, by Lemma \ref{f1}, we have $\sum_{F\in \mathcal{N}_2} \eta(F) \leq
\eta^+(G[U])\leq g_0$.
Therefore, $\sum_{F\in \mathcal{N}_2} |V(F)|$ is 
upper bounded by $g_0/\big(t_0-\frac1{\lfloor1/t_0\rfloor+1}\big).$
\end{proof}
 
\begin{claim}\label{cc4}
Every component in $\mathcal{N}_2\cup \mathcal{N}_3$ is a star. 
\end{claim}

\begin{proof}
Suppose to the contrary that $\mathcal{N}_2\cup \mathcal{N}_3$ 
contains some components that are not stars.
Then, let $G'$ be the graph obtained from $G$ by replacing
each component in $\mathcal{N}_2\cup \mathcal{N}_3$ with a star of the same order.
By Lemma \ref{l2},
the star is strictly walk-preferable to other trees with the same order.
By Claim \ref{cc3} and the definition of $\mathcal{N}_3$,
every component in $\mathcal{N}_2\cup \mathcal{N}_3$ has finite order.
Hence, by Corollary \ref{f0}, 
we have $\rho(G')>\rho(G)$. 

However, in view of (\ref{z00-1}), 
we know that $\eta(F)= t_0 |V(F)| - 1$ for each tree $F\in \mathcal{N}_2\cup \mathcal{N}_3$.
This implies that $\eta^+(G'[U])=\eta^+(G[U])\leq g_0$.
By Lemma \ref{f1}, $G'$ satisfies the property $P_{t,r}$,
which leads to a contradiction.
 \end{proof}
 
From (\ref{z00}), it follows that $\eta(F)\leq 0$ for every tree $F$ with $|V(F)|\leq\lfloor1/t_0\rfloor$.
Hence, no component in $\mathcal{N}_3$ is positive.
Furthermore, the following statement holds for $\mathcal{N}_3$. 

\begin{claim}\label{cc5}
There is at most one star $F$ in $\mathcal{N}_3$ with $|V(F)|\leq\lfloor1/t_0\rfloor-1$. 
\end{claim}

\begin{proof}
Suppose, for contradiction, that there are two stars $K_{1,a}$ and $K_{1,b}$ in $\mathcal{N}_3$ 
such that $0\leq b \leq a \leq \lfloor1/t_0\rfloor - 2$. 
If $b\geq1$, we construct a graph $G'$ from $G$ by performing an edge-shifting operation, 
i.e., replacing $K_{1,a} \cup K_{1,b}$ with $K_{1,a+1} \cup K_{1,b-1}$. 
In light of (\ref{z00}), neither $K_{1,a+1}$ nor $K_{1,b-1}$ are positive.
Thus, $\eta^+(G'[U])=\eta^+(G[U])\leq g_0$.
By Lemma \ref{f1}, $G'$ satisfies the property $P_{t,r}$.
Notice that $F_2(K_{1,a+1}\cup K_{1,b-1})>F_2(K_{1,a}\cup K_{1,b})$.
By Corollary \ref{f0}, we have $\rho(G')>\rho(G)$, 
which contradicts the maximality of $\rho(G)$. 

If $b = 0$, i.e., $K_{1,b}$ is an isolated vertex, 
we construct a new graph $G''$ from $G$ by adding an edge 
between the center of $K_{1,a}$ and the isolated vertex. 
It is clear that $\rho(G'')>\rho(G)$. 
Moreover, $\eta(K_{1,a+1})\leq0$,
and hence $\eta^+(G''[U])=\eta^+(G[U])\leq g_0$.
Therefore, $G''$ satisfies the property $P_{t,r}$,
which also leads to a contradiction.
\end{proof}

If $\mathcal{N}_2$ contains two stars of order at least $\lfloor1/t_0\rfloor+2$,
we can also perform an edge-shifting operation, i.e., 
replacing $K_{1,\lfloor1/t_0\rfloor + a}\cup K_{1, \lfloor1/t_0\rfloor+b}$ 
with $K_{1,\lfloor1/t_0\rfloor+a+1}\cup K_{1,\lfloor1/t_0\rfloor+b-1}$.
Using a similar argument as in Claim \ref{cc5}, we obtain the following statement.

\begin{claim}\label{cc6}
There is at most one star $K_{1,\lfloor\frac{1}{t_0}\rfloor+a}$ in $\mathcal{N}_2$ with $a \geq 1$.
\end{claim}

By Claims \ref{cc2} and \ref{cc3}, 
$\sum_{F\in \mathcal{N}_3}|V(F)|=\Theta(n)$.
Furthermore, by Claim \ref{cc5}, $\mathcal{N}_3$ contains $\Theta(n)$
copies of $K_{1,\lfloor1/{t_0}\rfloor-1}$.
Based on this, we have the following assertion.

\begin{claim}\label{cc7}
There is at most one component $F$ in $\mathcal{N}_1$.
If $|\mathcal{N}_1|=1$, then $F$ is a triangle. 
\end{claim}

\begin{proof}
Suppose to the contrary that there are two components $F_1, F_2 \in \mathcal{N}_1$. 
Then $e(F_i)\geq |V(F_i)|\geq3$ for $i\in\{1,2\}$. 
Let $\sum_{i=1}^2e(F_i)=\sum_{i=1}^2|V(F_i)|- 1 + a$.
By Claim \ref{cc2}, we obtain $\sum_{i=1}^2|V(F_i)|\leq g_0/t_0$,
and thus $a$ is also bounded by a constant.
Recall that $\mathcal{N}_3$ contains $\Theta(n)$ stars of order $\lfloor1/t_0\rfloor$. 
Now, consider the subgraph $F$ formed by the union of $F_1$, $F_2$, 
and $a$ stars $K_{1, \lfloor1/t_0\rfloor-1}$ from $\mathcal{N}_3$. 
It is clear that $e(F)=|V(F)|-1$, and $|V(F)|$ is finite.
In view of (\ref{z00}), we have
$\eta(K_{1,\lfloor1/t_0\rfloor-1})\leq 0$.
Recall that $t_0=1+\lfloor t\rfloor-t$.
Then 
\begin{align}\label{z00-3}
\eta^+(F)=\sum_{i=1}^2\eta(F_i)
=\sum_{i=1}^2e(F_i)+(\lfloor t\rfloor-t)\sum_{i=1}^2|V(F_i)|
=t_0\sum_{i=1}^2|V(F_i)|-1+a.
\end{align}

Now, we transform $F$ into a star $F' = K_{1, |V(F)| - 1}$
and construct a graph $G'$ from $G$ by replacing the edges of $F$ with the edges of $F'$. 
From (\ref{z00-1}),
we observe that $\eta^+(F')=\eta(F')=t_0|V(F')|-1$,
where $|V(F')|=\sum_{i=1}^2|V(F_i)|+a\lfloor{1}/{t_0}\rfloor$.
By comparing with (\ref{z00-3}), it is straightforward to verify that $\eta^+(F')\leq \eta^+(F)$.
Thus, $\eta^+(G'[U])\leq\eta^+(G[U])\leq g_0$, and hence $G'$ satisfies the property $P_{t,r}$.
Notice that $e(F')=e(F)>3$. 
By Lemma \ref{l2}, we see that $F_2(F')>F_2(F)$. 
Since $|V(F)|$ is finite, by Corollary \ref{f0}, we obtain $\rho(G')>\rho(G)$, 
which leads to a contradiction. 
Therefore, $|\mathcal{N}_1|\leq1$.

Next, assume that $|\mathcal{N}_1|=\{F_1\}$. We will demonstrate that $F_1$ is a triangle.
The proof follows similarly to the previous one.
If $e(F_1)\geq 4$, we denote $e(F_1)=|V(F_1)|-1+a$ and $F=F_1\cup aK_{1,\lfloor 1/t_0\rfloor-1}$.
Let $G'$ be a graph obtained from $G$ by removing the edges of $F$ and adding the edges of $F'$,
where $F'=K_{1,|V(F)|-1}$.
We also have $\eta^+(F')\leq\eta^+(F),$ and hence
$\eta^+(G'[U])\leq\eta^+(G[U])\leq g_0$, 
which implies that $G'$ satisfies the property $P_{t,r}$.
By a similar argument as above,
we can conclude that $\rho(G')>\rho(G)$,
a contradiction.
Therefore, $e(F_1)=3$, and so $F_1$ is a triangle. 
This completes the proof.
\end{proof}

With the constraints on $\mathcal{N}_1,\mathcal{N}_2$ and $\mathcal{N}_3$ established as outlined above,
we can now relate these findings to the overall structure of
$G[U]$. Specifically,
we will investigate the implications of having a single component in $\mathcal{N}_1$.

\begin{claim}\label{cc8}
 The following statements are equivalent:\\
(i) $|\mathcal{N}_1| = 1$;\\
(ii) $3t_0 \leq g_0 < 8t_0 - 4$;\\
(iii) $G[U] \cong K_3 \cup (n - \lfloor t \rfloor - 3)K_1$.
\end{claim}

\begin{proof}
We prove the equivalence of the three statements 
by showing the implications 
$(i) \Rightarrow (ii)$, $(ii) \Rightarrow (iii)$, and $(iii) \Rightarrow (i)$.

$(i) \Rightarrow (ii)$:
Assume that $|\mathcal{N}_1| = 1$. By Claim \ref{cc7}, 
the unique component in $\mathcal{N}_1$ is a triangle $K_3$. 
Recall that $t_0 =\lfloor t \rfloor+ 1 - t> 0$.
Then $\eta(K_3) = e(K_3) + (\lfloor t \rfloor - t)|V(K_3)| = 3t_0$. 
Since $\eta(K_3) \leq \eta^+(G[U]) \leq g_0$, we have $g_0 \geq 3t_0$.

Next, we show that $g_0 < 8t_0 - 4$. Suppose, for contradiction, that $g_0 \geq 8t_0 - 4$. 
Since $\mathcal{N}_1$ contains a triangle, we have $e(G[U])\geq4$. 
Otherwise, $\mathcal{N}_2 \cup \mathcal{N}_3$ would consist only of isolated vertices. 
Recall that $\mathcal{N}_3$ contains $\Theta(n)$ components.
We can construct a new graph $G'$ from $G$ by replacing $K_3 \cup 5K_1$ with $4K_2$. 
We observe that $\eta^+(G'[U]) \leq \eta(4K_2) = 8t_0 - 4 \leq g_0$.
Thus, by Lemma \ref{f1}, $G'$ satisfies the property $P_{t,r}$. 
However, since $4K_2$ is walk-preferable to $K_3 \cup 5K_1$, 
Corollary \ref{f0} implies that $\rho(G') > \rho(G)$, a contradiction. 
This contradiction establishes the existence of a star 
$K_{1,a}$ in $\mathcal{N}_2 \cup \mathcal{N}_3$ with $a \geq 1$.

Now, we construct another graph $G''$ from $G$ by replacing $K_3 \cup K_{1,a}$ with $K_{1,3+a}$. 
By Lemma \ref{l2} and Corollary \ref{f0}, we have $\rho(G'') > \rho(G)$. Since
$K_{1,3+a}$ and $K_3 \cup K_{1,a}$ have the same number of vertices and the same number of edges,
we have $\eta(K_{1,3+a}) = \eta(K_3 \cup K_{1,a})$,
and hence, $\eta^+(G''[U])=\eta^+(G[U]) \leq g_0.$
Thus, $G''$ also satisfies the property $P_{t,r}$, 
contradicting the maximality of $\rho(G)$. 
Therefore, $g_0 < 8t_0 - 4$.

$(ii) \Rightarrow (iii)$:
Assume that $3t_0 \leq g_0 < 8t_0 - 4$. 
Then $t_0>\frac45$.
If $e(G[U])\leq2$, we
modify the graph $G$ by removing the edges within $U$ and then
embedding a triangle in $U$.
By a similar argument as above, we get a contradiction.
Thus, $e(G[U])\geq3$. 

Let $F$ be the subgraph induced by all the edges in $G[U]$.
Clearly, $|V(F)|\leq 2e(F)$, and thus,
$\eta(F)=e(F)-(1- t_0)|V(F)|\geq e(F)(2t_0-1).$
If $e(F)\geq4$, then 
$\eta(F)\geq 8t_0 - 4$ and so $\eta^+(G[U])=\eta(F)>g_0$. 
This is a contradiction. 
Thus, we have $e(F)=3$.
Furthermore, 
by Lemma \ref{l2}, $F$ must be either $K_{1,3}$ or $K_3$. 
Since $W_3(K_3) = 24 > W_3(K_{1,3}) = 18$, $K_3$ is walk-preferable to $K_{1,3}$. 
Therefore, $F \cong K_3$, and so $G[U] \cong K_3 \cup (n - \lfloor t \rfloor - 3)K_1$.

$(iii) \Rightarrow (i)$: If $G[U] \cong K_3 \cup (n - \lfloor t \rfloor - 3)K_1$, 
then $G[U]$ contains exactly one non-tree component $K_3$. Hence, $|\mathcal{N}_1| = 1$.

This completes the proof.
\end{proof}

Building on the preceding analysis and Claim \ref{cc8}, we can now establish the following theorem,
which uniquely determines the extremal graph $G$ 
when $3t_0 \leq g_0 < 8t_0 - 4$.

\begin{theorem}\label{t4}
Let $G$ be a graph that attains the maximum spectral radius over all
$n$-vertex graphs satisfying the property $P_{t,r}$. 
If $3t_0 \leq g_0 < 8t_0 - 4$, then we have
\[
G \cong K_{\lfloor t \rfloor} \nabla \left(K_3 \cup (n - \lfloor t \rfloor - 3)K_1\right).
\]
\end{theorem}

In the following, we need only consider the case 
where $g_0$ does not satisfy the inequality 
$3t_0 \leq g_0 < 8t_0 - 4 $.
By Claims \ref{cc7} and \ref{cc8}, $\mathcal{N}_1$ must be empty. 
Thus, $G[U]$ is a forest. Furthermore, 
$G[U]$ consists solely of stars by Claim \ref{cc4}.

According to the definition of $\mathcal{N}_2$, 
the stars within $\mathcal{N}_2$ are of the form $K_{1,\lfloor1/t_0\rfloor + a},$ 
where $a \geq 0$. By Claim \ref{cc6}, there are only two distinct types of stars in $\mathcal{N}_2$. 
We refer to the star as a {\em big star} when $a = 0$, 
and as a {\em maximal star} for cases where $a \geq 1$.
For $\mathcal{N}_3$, the stars are of the form $K_{1, \gamma-1}$, 
where $1 \leq \gamma \leq \lfloor1/t_0\rfloor $. 
These stars are referred to as {\em small stars} 
when $\gamma = \lfloor1/t_0\rfloor$, 
and as {\em minimal stars} when $1 \leq \gamma <\lfloor1/t_0\rfloor$. 
By Claims \ref{cc5} and \ref{cc6},
there exists at most one minimal (or maximal) star. 
In particular, when $\gamma=1$, a minimal star becomes an isolated vertex.

To determine the number of edges in $G[U]$, 
we introduce the following notations: 
Let $b$ denote the number of big stars (i.e., $K_{1, \lfloor1/t_0\rfloor}$) in $G[U]$. 
For convenience, when considering the maximal star $K_{1,\lfloor1/t_0\rfloor+ a}$, 
allow $a=0$. If $a = 0$, then the number of big stars is $b + 1$; 
otherwise, there is exactly one maximal star of the form 
$K_{1,\lfloor1/t_0\rfloor+ a}$ with $a \geq 1$.

Define $\phi_0 = \eta(K_{1, \lfloor1/t_0\rfloor}) = (1 + \lfloor1/t_0\rfloor)t_0 - 1$, 
noting that $0 < \phi_0 \leq t_0$. Express $g_0$ as
\begin{equation}\label{e7-222}
g_0 = \alpha_0 \phi_0 + \beta_0,
\end{equation}
where $\alpha_0$ is a nonnegative integer and $0 \leq \beta_0 < \phi_0$.

Additionally, let $n - \lfloor t \rfloor - \alpha_0 
= n_0 \lfloor1/t_0\rfloor+ \gamma_0$, 
where $n_0$ is an integer sufficiently large, and $1 \leq \gamma_0 \leq \lfloor1/t_0\rfloor$.
Using these definitions, we now proceed to determine $e(G[U])$.

\begin{claim}\label{cc9}
$e(G[U]) = n - \lfloor t \rfloor - (n_0 + 1)$.
\end{claim}

\begin{proof}
By Claim \ref{cc3}, the stars within $\mathcal{N}_2$ are positive. 
Moreover, recall that $\mathcal{N}_3$ contains no positive stars.
Then,  $\eta^+(G[U])=\eta(K_{1,\lfloor1/t_0\rfloor+a})+b\eta(K_{1,\lfloor1/t_0\rfloor}).$
It follows that
\begin{equation}\label{e7-111}
\eta^+(G[U]) = a t_0 + (b + 1) \phi_0 = (a + b + 1) \phi_0
+ a (t_0-\phi_0) \geq (a + b + 1) \phi_0.
\end{equation}
Note that $\eta^+(G[U])\leq g_0$.
If $a + b + 1 \geq \alpha_0 + 1$, then $g_0 \geq (\alpha_0+1)\phi_0 > \alpha_0 \phi_0 + \beta_0 = g_0$,
which is a contradiction. Therefore, 
\begin{equation}\label{e7-1111}
a + b + 1 \leq \alpha_0.
\end{equation}

Since $G[U]$ is a forest, and based on the classification of $\mathcal{N}_2\cup\mathcal{N}_3$,
the number of components in $G[U]$ is
$\big\lceil \frac{n - \lfloor t \rfloor - (a + b + 1)}{\lfloor1/t_0\rfloor} \big\rceil$.
Thus, the number of edges in $G[U]$ is
\begin{equation}\label{e7}
e(G[U]) = n - \lfloor t \rfloor - \Big\lceil \frac{n - \lfloor t \rfloor -
(a + b + 1)}{\lfloor1/t_0\rfloor} \Big\rceil .
\end{equation}
Combining \eqref{e7-1111} and \eqref{e7}, we obtain that
\[
e(G[U]) \leq n - \lfloor t \rfloor - \Big\lceil \frac{n - \lfloor t \rfloor
- \alpha_0}{\lfloor1/t_0\rfloor} \Big\rceil = n - \lfloor t \rfloor - (n_0 + 1).
\]

Next, we show that
$e(G[U]) \geq n - \lfloor t \rfloor - (n_0 + 1)$.
Suppose, to the contrary, that this is not the case.
Let $F$ be the union of the maximal star, the minimal star, all big stars, and some (finite) small stars.
By edge-shifting operations, we can transform $F$ into a forest $F'$,
which is the union of $\alpha_0$ big stars and a new minimal star. Then construct a graph $G'$ from $G$ by 
replacing the edges of $F$ with the edges of $F'$. 
It is clear that $\eta^+(G'[U]) = \alpha_0 \phi_0 \leq g_0$,
which implies that $G'$ satisfies the property $P_{t,r}$. Moreover,
\[
e(G'[U]) = n - \lfloor t \rfloor - \Big\lceil \frac{n - \lfloor t \rfloor
- \alpha_0}{\lfloor1/t_0\rfloor} \Big\rceil =n - \lfloor t \rfloor - (n_0 + 1) \geq e(G[U]) + 1.
\]
By Corollary \ref{f0}, we have $\rho(G') > \rho(G)$, a contradiction. This completes the proof.
\end{proof}

Next, we determine the values of $a$ and $b$. To achieve this, we shall introduce a new parameter.
Denote $c=\alpha_0-(a+b+1)$. Based on \eqref{e7-1111}, it follows that $c \geq 0$.
Recall that $n - \lfloor t \rfloor = \alpha_0 + \gamma_0+ n_0 \lfloor1/t_0\rfloor$. 
Combining this with equality \eqref{e7} and Claim \ref{cc9}, 
we have $ \alpha_0+\gamma_0-(a + b + 1) \leq \lfloor1/t_0\rfloor$.
It follows that
$c \leq \lfloor1/t_0\rfloor - \gamma_0$.

We can observe that $G$ is uniquely determined by $a$ and $b$.
Therefore, we may denote $G = G_{(a,b,c)}$, and let $\eta^+(G_{(a,b,c)}[U]) = \eta^+_{(a,b,c)}$. 
Recall that $\eta^+_{(a,b,c)}=\eta^+(G[U])\leq g_0$.
Furthermore,
based on \eqref{e7-111}, we obtain that
\begin{equation}\label{e7-11111}
\eta^+_{(a,b,c)}=at_0+(b+1)\phi_0
=a t_0+(\alpha_0-a-c)\phi_0 \leq g_0 = \alpha_0 \phi_0 + \beta_0.
\end{equation}
From the previous discuss, we get the following inequalities:
\begin{eqnarray}
  a t_0 - (a + c) \phi_0 &\leq& \beta_0, \label{e8} \\
  0 \leq c &\leq& \lfloor1/t_0\rfloor - \gamma_0, \label{e9} \\
  a + c &\leq& \alpha_0 - 1. \label{e10}
\end{eqnarray}

In view of \eqref{e7-11111}, we have
$\eta^+_{(0, \alpha_0-1, 0)} = \alpha_0 \phi_0 = g_0 - \beta_0,$
and by Lemma \ref{f1}, we see that $G_{(0,\alpha_0-1, 0)}$ satisfies the property $P_{t,r}$.
Observe that $G_{(0, \alpha_0-1, 0)}$ contains $\alpha_0$ components of the form $K_{1,\lfloor1/t_0\rfloor}$.
Additionally, recall that $n - \lfloor t \rfloor - n_0 \lfloor1/t_0\rfloor=\alpha_0 + \gamma_0$.
This further implies that $G_{(0, \alpha_0-1, 0)}$
contains a minimal component $K_{1,\gamma_0-1}$ and $n_0-\alpha_0$ components of the form $K_{1,\lfloor1/t_0\rfloor-1}$.
Recall that $\alpha_0=a+b+c+1$. It is easy to see that
$G_{(a,b,c)}$ can be obtained from $G_{(0, \alpha_0-1, 0)}$
by the following two types of edge-shifting operations:
\vspace{1.5mm}
 
 (i)  $F_1 := (a+1)K_{1,\lfloor1/t_0\rfloor}
  \implies F'_1 := K_{1,\lfloor1/t_0\rfloor+a} \cup aK_{1,\lfloor1/t_0\rfloor-1}$;
  
 (ii)  $F_2 := cK_{1,\lfloor1/t_0\rfloor}
  \cup K_{1,\gamma_0-1} \implies F'_2 := cK_{1,\lfloor1/t_0\rfloor - 1} \cup K_{1,c+\gamma_0-1}$.

\vspace{1.5mm}
We now define $f_k(a,c) = W_k(F'_1 \cup F'_2) - W_k(F_1 \cup F_2)$.
Then, it is clear that $f_1(a,c)=0$ since $e(F'_1\cup F'_2)=e(F_1\cup F_2)$. 
Note that $W_2(F)=\sum_{v \in V(F)}\big(d_F(v)\big)^2$. We further have
\begin{equation*}
\begin{aligned}
f_2(a,c) &= (a + \lfloor1/t_0\rfloor)^2 + (c + \gamma_0 - 1)^2 + (a + c)(\lfloor1/t_0\rfloor - 1)^2 \\
&\quad - (a + c + 1)\lfloor1/t_0\rfloor^2 - (\gamma_0 - 1)^2 \\
&= a^2 + 2a\lfloor1/t_0\rfloor + c^2 + 2c(\gamma_0 - 1) - 2(a + c)\lfloor1/t_0\rfloor + (a + c) \\
&= a^2 + a + c^2 - 2c\big(\lfloor1/t_0\rfloor - \gamma_0+1/2\big).
\end{aligned}
\end{equation*}

Subject to the property $P_{t,r}$, 
the extremal graph $G$ is obtained precisely when the function $f_2(a,c)$ reaches its maximum.
Therefore, our goal is to determine the values of $a$ and $c$ that maximize $f_2(a, c)$,
under the condition that $\eta^+_{(a,b,c)}\leq g_0$.

\begin{claim}
$a=\min\big\{\big\lfloor\frac{\beta_0}{t_0-\phi_0}\big\rfloor,\alpha_0-1\big\}$,
$b=\alpha_0-a-1$ and $c=0$.
\end{claim}

\begin{proof}
First, observe that $f_2(a, c)$ is a quadratic function with respect to $a$ and $c$, 
and it represents an upward-opening paraboloid of revolution.
The symmetry plane of $f(a, c)$ with respect to $c$ is $c = \big\lfloor1/t_0\big\rfloor - \gamma_0 + \frac{1}{2}$.
Form \eqref{e9}, we know that $0 \leq c\leq \lfloor1/t_0\rfloor - \gamma_0$.
Then, the extremum of $f_2(a, c)$ must occur when $c = 0$.

Next, from \eqref{e8} and \eqref{e10}, we can derive that
$a \leq \frac{\beta_0}{t_0 - \phi_0}$ and $a\leq\alpha_0-1$ when $c=0$.
Since $a$ is a nonnegative integer and the symmetry plane of $f_2(a,c)$ with respect to $a$ 
is $a =-\frac{1}{2}$, it follows that
$a = \min\big\{\big\lfloor\frac{\beta_0}{t_0 - \phi_0}\big\rfloor, \alpha_0 - 1\big\},$
where $a = \alpha_0 - 1$ when $t_0 = \phi_0$.
Recall that $c=\alpha_0-(a+b+1)$. 
We further have $b=\alpha_0 - a- 1$.

Now it suffices to verify that $\eta^+_{(a,b,c)}\leq g_0$ for the values of $a,b,c$ obtained above.
In view of \eqref{e7-11111}, we have
$$\eta^+_{(a,b,c)}=at_0+(b+1)\phi_0=
\alpha_0\phi_0+a(t_0-\phi_0)\leq\alpha_0 \phi_0 + \beta_0= g_0.$$
This completes the proof.
\end{proof}

Recall that $n-\lfloor t \rfloor- \alpha_0 = n_0\lfloor1/t_0\rfloor+\gamma_0$,
where $1\leq \gamma_0\leq \lfloor1/t_0\rfloor$. Based on the above arguments,
the following theorem holds.

\begin{theorem}\label{t5}
Let $G$ be a graph that attains the maximum spectral radius over all
$n$-vertex graphs satisfying the property $P_{t,r}$. 
If $g_0$ does not satisfies $3t_0 \leq g_0 < 8t_0 - 4$,
then $G\cong K_{\lfloor t \rfloor}\nabla H$,
where $H $ consists of stars as follows:

(i) one maximal star $K_{1, \lfloor1/t_0\rfloor+a}$,
  where $a= \min\big\{\big\lfloor\frac{\beta_0}{t_0-\phi_0}\big\rfloor,\alpha_0-1\big\}$;
  
(ii) $b$ big stars $K_{1,\lfloor1/t_0\rfloor}$, where $b=\alpha_0-a-1$;

(iii) $\ell$ small stars $K_{1,\lfloor1/t_0\rfloor-1}$, where $\ell=n_0-b-1$;

(iv) one minimal star $K_{1,\gamma-1}$, where $\gamma=\gamma_0 $.
\end{theorem}

Note that when $a=0$, it indicates that there is no maximal star in $G[U]$.
Consequently, the number of big stars in
$G[U]$ simplifies to $\alpha_0$. In the case $\gamma_0=\lfloor1/t_0\rfloor$,
it implies that there is no minimal star in $G[U]$. 
The absence of a minimal star changes the composition of small stars.
Specifically, the number of small stars becomes
$n_0-b$.

In summary, by using Theorems \ref{t4} and \ref{t5}, 
we have successfully completed the proof of Theorem \ref{t3}. 
To conclude this section, we will present a corollary.

If $t$ is a positive integer, then $t_0=\lfloor t\rfloor+1-t=1$ and
$\phi_0=(1 + \lfloor1/t_0\rfloor)t_0-1=1$.
Based on \eqref{e7-222}, we further have $\alpha_0=\lfloor g_0\rfloor.$
It follows that $a = \alpha_0 - 1=\lfloor g_0\rfloor-1$ and $b = 0$. 
By applying Theorems \ref{t4} and \ref{t5}, we obtain the following corollary.

\begin{cor}\label{integer}
Let $G$ be a graph that attains the maximum spectral radius over all
$n$-vertex graphs satisfying the property $P_{t,r}$. 
If $t$ is a positive integer, then 
$G\cong K_t\nabla \big(K_3 \cup (n-t- 3)K_1\big)$ for $3\leq g_0<4$,
or $G \cong K_t\nabla \big(K_{1,\lfloor g_0\rfloor} \cup (n - t- \lfloor g_0\rfloor- 1)K_1\big)$
otherwise.
\end{cor}

For $t$ being a positive integer, Guiduli proved the Theorem \ref{t0}, 
and he also conjectured that the spectral extremal graph with respect to Theorem \ref{t0} is 
$K_t\nabla F$ where the edges within $F$ form a star. Corollary \ref{integer} resolves this conjecture, 
and provides a complete characterization of the extremal graphs.  

\section{Concluding remarks}

We may also consider a generalization of the Hereditarily Bounded Property,
where the hereditary bound on the number of edges is not restricted to linear constraints.
In his thesis, Guiduli posed the following problem.
\begin{prob}
Given a graph $G$, if $|E(H)|\leq c\cdot|V(H)|^2$ holds for all $H\subseteq G$,
does it imply that $\rho(G)\leq 2c\cdot|V(G)|$? 
What can be said if the exponent of $2$ is replaced by some other constant less than $2$?
\end{prob}

This issue has also been highlighted by Liu and Ning in their survey of unsolved problems in spectral graph theory \cite{liu23}.
It is worth noting that $\rho(G)\leq\sqrt{2c}n$ serves as a trivial upper bound, and that $2c$ would be best possible,
as seen by Wilf's inequality \cite{Wilf86}, i.e., $\rho(G)\leq (1-\frac{1}{\omega})n$, where $\omega$ denotes the clique number of $G$. If this question were true,
then the spectral Erd\H{o}s-Stone-Simonovits theorem \cite{Nikiforov09}
would follow from the Erd\H{o}s-Stone-Simonovits theorem \cite{Erdos46, Erdos66}
and Wilf's theorem would follow from Tur\'{a}n's theorem \cite{Turan61}.

\vspace{-0.2cm}

\end{document}